\documentclass{amsart}
\usepackage{graphicx}
\newtheorem{theorem}{Theorem}[section]
\newtheorem{lemma}[theorem]{Lemma}
\newtheorem{corollary}[theorem]{Corollary}

\theoremstyle{definition}

\theoremstyle{remark}

\numberwithin{equation}{section}

\begin{document}

\title[Closed incompressible surfaces in the complements of positive knots]{Closed incompressible surfaces in the complements of positive knots}

\author{Makoto Ozawa}
\address{Department of Mathematics, School of Education, Waseda University, Nishiwaseda 1-6-1, Shinjuku-ku, Tokyo 169-8050, Japan}
\email{ozawa@musubime.com}
\thanks{The author was supported in part by Fellowship of the Japan Society for the Promotion of Science for Japanese Junior Scientists.}

\date{}


\begin{abstract}
We show that any closed incompressible surface in the complement of a positive knot is algebraically non-split from the knot, positive knots cannot bound non-free incompressible Seifert surfaces and that the splitability and the primeness of positive knots and links can be seen from their positive diagrams.
\end{abstract}

\maketitle

\section{Introduction}
A knot $K$ in the 3-sphere $S^3$ is called {\it positive} if it has an oriented diagram all crossings of which are positive crossings.
For a closed surface $F$ in $S^3-K$, we define the {\it order} $o(F;K)$ of $F$ for $K$ as follows (\cite{O}).
Let $i:F\to S^3-K$ be the inclusion map and let $i_*:H_1(F)\to H_1(S^3-K)$ be the induced homomorphism.
Since $Im(i_*)$ is a subgroup of $H_1(S^3-K)=\Bbb{Z} \langle$meridian$\rangle$, there is an integer $m$ such that $Im (i_*)= m \Bbb{Z}$. Then we define $o(S;K)=m$. 

The positive knot complements have the following special properties.

\begin{theorem}
Any closed incompressible surface in a positive knot complement has non-zero order.
\end{theorem}

A Seifert surface $F$ for a knot is said to be {\it free} if $\pi _1(S^3-F)$ is a free group.
In \cite[Theorem 1.1]{O}, it is shown that a knot bounds a non-free incompressible Seifert surface if and only if there exists a closed incompressible surface in the knot complement whose order is equal to zero. Therefore, Theorem 1.1 gives us the next corollary.

\begin{corollary}
Positive knots cannot bound non-free incompressible Seifert surfaces.
\end{corollary}

Although positive links which have connected positive diagrams are non-split because they have positive linking numbers, we can give another geometrical proof of this fact.

\begin{theorem}
Positive links are non-split if their positive diagrams are connected.
\end{theorem}

Positive diagrams of positive knots or links also tell us their primeness.
We say that a knot or link diagram $\tilde{K}$ on the 2-sphere $S$ is {\it prime} if for any loop $l$ in $S$ intersecting $\tilde{K}$ in 2 points, $l$ bounds a disk intersecting $\tilde{K}$ in an arc.

\begin{theorem}
Non-trivial positive knots or links are prime if their positive diagrams are connected and prime.
\end{theorem}

Theorem 1.4 widens the result of Cromwell (\cite[1.2 Theorem]{C}).

\section{Proof of Theorem 1.1 and 1.3}
Theorem 1.1 and 1.3 follow the next Theorem.


\begin{theorem}
Let $K$ be a positive knot or link in the 3-sphere $S^3$ and  $F$ a closed incompressible surface in the complement of $K$.
Then one of the following conclusions (1) and (2) holds.

\begin{enumerate}
	\item[(1)] There exists a loop $l$ in $F$ such that $lk(l,K)\ne 0$.
	\item[(2)] $F$ is a splitting sphere for $K$, and any positive diagram of $K$ is disconnected.
\end{enumerate}
\end{theorem}

Henceforth, we shall prove Theorem 2.1.

Let $S$ be a 2-sphere in $S^3$ and $p:S^3-\{ 2\ points\} \cong S\times R \to S$ a projection.
Put $K$ so that $p(K)$ is a positive diagram.
As usual way, we express $K$ in a bridge presentation.
Thus we have the following data (see Figure 1).

\begin{itemize}
	\item $S^3=B^+\cup _S B^-$ ($S$ decomposes $S^3$ into two 3-balls)
	\item $K=K^+\cup _S K^-$, where $K^{\pm}\subset B^{\pm}$ ($S$ cuts $K$ into over bridges and under bridges)
	\item $K^{\pm}=K^{\pm} _1 \cup K^{\pm} _2 \cup \ldots K^{\pm} _n$ ($K$ is presented as $n$ over bridges and $n$ under bridges)
	\item $D^{\pm}=D^{\pm}_1 \cup D^{\pm}_2 \cup \ldots D^{\pm}_n$ (each $K^{\pm}_i \cup p(K^{\pm}_i)$ bounds a disk $D^{\pm}_i$ such that $p(D^{\pm}_i)=p(K^{\pm}_i)$)
\end{itemize}

\begin{figure}[htbp]
	\begin{center}
		\includegraphics[trim=0mm 0mm 0mm 0mm, width=.6\linewidth]{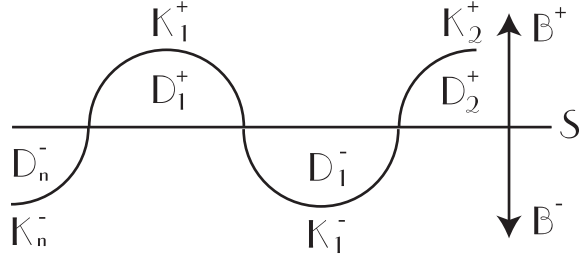}
	\end{center}
	\caption{View from level surface}
\end{figure}

We take $n$ minimal over all bridge presentations of $p(K)$.

\begin{lemma}
We may assume that $F\cap D^-=\emptyset$, $F\cap B^-$ consists of disks, $F\cap D^+$ consists of arcs, and any component of $F\cap B^+-D^+$ is a disk.
\end{lemma}

\begin{proof}
This can be done by cutting and pasting $F$ along some disks.
Note that such operations do not have any effect on the conditions (1) and (2) if we take a suitable choice of $F$.
\end{proof}

We take $|F\cap B^-|$ and $|F\cap D^+|$ minimal.
Note that $|F\cap B^-|\ne 0$ because $F$ is incompressible in $S^3-K$.
If $|F\cap B^-|=1$ and $|F\cap D^+|=0$, then we have the conclusion (2).

Hereafter, we suppose that $|F\cap B^-|\ge 1$ and $|F\cap D^+|\ge 1$.

Then we obtain a connected graph $G$ in $F$ by regarding $F\cap B^-$ and $F\cap D^+$ as vertices and edges respectively.
Note that every vertex has a positive even valency by the construction.

An arc $\alpha _j$ of $F\cap D^+_i$ divides $D^+_i$ into two disks $\delta _j$ and $\delta _j'$, where $\delta _j'$ contains $K^+_i$.
Put $\beta _j=\delta _j\cap S$.
We may assume that $p(\alpha _j)=p(\delta _j)=\beta _j$ for all $\alpha _j$.
We assign an orientaion endowed from $K_i$ to $\alpha _j$ and $\beta _j$ naturally (see Figure 2).

\begin{figure}[htbp]
	\begin{center}
		\includegraphics[trim=0mm 0mm 0mm 0mm, width=.5\linewidth]{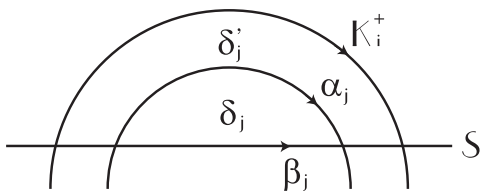}
	\end{center}
	\caption{$\alpha _j$ and $\beta _j$ have the orientaions}
\end{figure}

\begin{lemma}
For any arc $\alpha _j$ of $F\cap D^+_i$, $\beta _j\cap p(K^-)\ne \emptyset$.
\end{lemma}

\begin{proof}
Suppose that there exists an arc $\alpha _j$ of $F\cap D^+_i$ such that $\beta _j\cap p(K^-)= \emptyset$.
By exchanging $\alpha _j$ if necessary, we may assume that $\alpha _j$ is outermost in $D^+_i$, that is,  $int\delta _j\cap F=\emptyset$.
If $\alpha _j$ connects different vertices, then a $\partial$-compression of $F$ along $\delta _j$ reduces $|F\cap B^-|$.
Otherwise, $\alpha _j$ incidents a single vertex, say $D^-_k$.
We perform a $\partial$-compression of $F$ along $\delta _j$, and obtain an annulus $A$ consisting of the disk $D^-_k$ and the resultant band $b$.
Since we chose an outermost arc $\alpha _j$ and $\beta _j\cap p(K^-)= \emptyset$, there exists a compressing disk for $A$ in $B^--K^-$.
By retaking $F$ along the compressing disk, we can reduce $|F\cap D^+|$.
\end{proof}

Now we pay attention to a face $f$ of $G$ in $F$.
The `cycle' $\partial f$ consists of edges and `corners' as subarcs in $\partial (F\cap B^-)$.
The edges have orientations as previously mentioned.

\begin{lemma}
For any face $f$, the cycle $\partial f$ can not be oriented.
\end{lemma}

\begin{proof}
Suppose that there is a face $f$ such that $\partial f$ can be oriented.
Then, since no corner of $\partial f$ intersects $p(K)$, and by Lemma 2.3, $p(\partial f)$ has non-zero intersection number with $p(k^-)$ on $S$ as illustrated in Figure 3.
This is a contradiction.
\end{proof}

\begin{figure}[htbp]
	\begin{center}
		\includegraphics[trim=0mm 0mm 0mm 0mm, width=.4\linewidth]{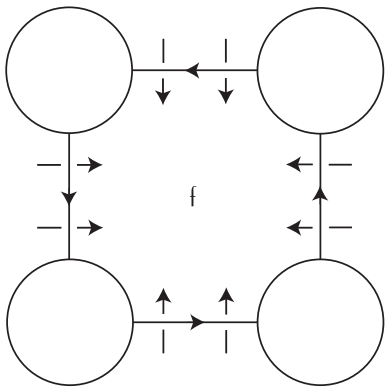}
	\end{center}
	\caption{$p(\partial f)$ has non-zero intersection number}
\end{figure}


For each face $f$ of $G$ and any point in the interior of any edge of $\partial f$, we can find an arc $\gamma$ on $f$ satisfying the following property.

\vspace{3mm}

(*) $\gamma$ connects two edges of $\partial f$ whose orientations are defferent in $\partial f$.
\vspace{3mm}


\begin{figure}[htbp]
	\begin{center}
		\includegraphics[trim=0mm 0mm 0mm 0mm, width=.4\linewidth]{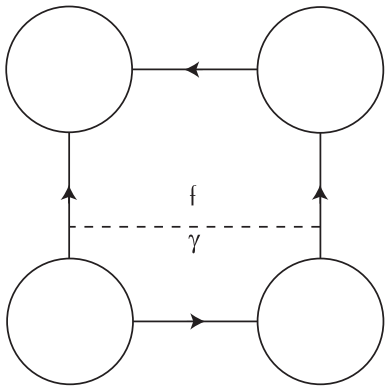}
	\end{center}
	\caption{$\gamma$ with the property (*)}
\end{figure}

Lemma 2.4 assures the existence of such an arc $\gamma$.

To find a loop $l$ on $F$ with $lk(l,K)\ne 0$, we depart a point in the interior of any edge of $G$, trace arcs with the property (*), and will arrive at the face on which we have walked.
Connecting these arcs, we will obtain an oriented loop $l$ in $F\cap B^+$ with a suitable orientaion such that $l$ has a positive intersection number with edges of $G$ on $F$.
Thus we got an oriented loop $l$ in $F$ which has non-zero linking number with $K$.
Since any loop in a splitting sphere is contractible in $S^3-K$, we have the conclusion (1).

This completes the proof of Theorem 2.1.

\section{Proof of Theorem 1.4}
Let $K$ be a positive knot or link in $S^3$ and $F$ be a decomposing sphere for $K$.
We put $K$ and $F$
as the proof of Theorem 1.1 except that two points $p_1$ and $p_2$ of $F\cap K$ are in $int B^+$ or $int B^-$.
Note that $p_1$ and $p_2$ can not be the ends of a single arc of $F\cap D^{\pm}$ because the tangle $(B^{\pm}, K^{\pm})$ is trivial and $F$ is a decomposing sphere.
Hence, there are two arcs $e_1$ and $e_2$ of $F\cap D^{\pm}$ whose ends contain $p_1$ and $p_2$ respectively.
We deform $F$ by an isotopy relative to $K$ so that $p(e_i)=p(p_i)$ $(i=1,2)$.
We take the number of bridges $n$ minimal.


\begin{lemma}
We may assume that $F\cap D^-\subset e_1\cup e_2$, $F\cap B^-$ consists of disks, $F\cap D^+$ consists of arcs, and any component of $F\cap B^+-D^+$ is a disk.
\end{lemma}

\begin{proof}
This can be done by an isotopy of $F$ since Theorem 1.3 assures us that $S^3-K$ is irreducible.
\end{proof}

We take $|F\cap B^-|$ and $|(F\cap D^+)-(e_1\cup e_2)|$ minimal.
Then we obtain a connected graph $G$ in $F$ by regarding $F\cap B^-$ and $(F\cap D^+)-(e_1\cup e_2)$ as vertices and edges respectively.
Corners of each face of $G$ may contain two points $\partial e_1-p_1$ and $\partial e_2-p_2$.
Note that $|F\cap B^-|\ne 0$, otherwise $F$ is not a decomposing sphere since $(B^{\pm},K^{\pm})$ is a trivial tangle.
If $|F\cap B^-|=1$ and $F\cap D^+\subset e_1\cup e_2$, then $F\cap S$ gives a desired loop since $p(e_i)=p(p_i)$ $(i=1,2)$.

\begin{lemma}
For any arc $\alpha _j$ of $(F\cap D^+)-(e_1\cup e_2)$, $\beta _j\cap p(K^-)\ne \emptyset$.
\end{lemma}

\begin{proof}
This can be done by the same argument to Lemma 2.3.
\end{proof}

Hereafter, we assume that $\tilde{K}$ is prime.

\begin{lemma}
There is no vertex of $G$ with valency 1.
\end{lemma}

\begin{proof}
Suppose that there is a vertex $V$ with valency 1.
Then only one edge $\alpha$ incident to $V$, and hence exactly one of $e_1$ and $e_2$ is attached to $V$ or contained in $V$.
Thus $\partial V$ intersects $\tilde{K}$ in two points.
Since $\tilde{K}$ is prime, $\partial V$ bounds a disk $E$ in $S$ which intersects $p(K)$ in an unknotted arc.
In the former case, $p(K)\cap E$ lies under a subarc of $K^+$ by the minimality of the number of bridges $n$.
Then by an isotopy of $F$ along the 3-ball which is bounded by $V\cup E$, we can reduce $|F\cap B^-|$. See Figure 5.
In the later case, $E$ intersects $K$ in one point, and $V\cup E$ bounds a pair of a 3-ball and an unknotted subarc of $K^-$ by the minimality of $n$.
Then an isotopy of $F$ along the pair can reduce $|F\cap B^-|$. See Figure 6.
\end{proof}

\begin{figure}[htbp]
	\begin{center}
		\includegraphics[trim=0mm 0mm 0mm 0mm, width=.7\linewidth]{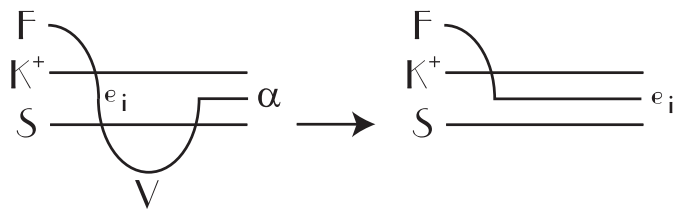}
	\end{center}
	\caption{Isotopy of $F$ along the 3-ball}
\end{figure}

\begin{figure}[htbp]
	\begin{center}
		\includegraphics[trim=0mm 0mm 0mm 0mm, width=.7\linewidth]{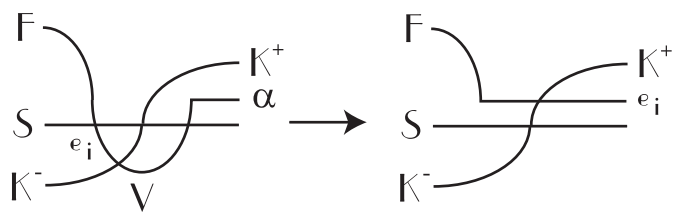}
	\end{center}
	\caption{Isotopy of $F$ along the pair}
\end{figure}

\begin{lemma}
There is no face $f$ of $G$ in $F$ such that $\partial f$ is a loop of $G$.
\end{lemma}

\begin{proof}
Suppose there exists a face $f$ as Lemma 3.4.
Then $\partial f$ consists of an edge $\alpha$ of $G$ and a subarc $\gamma$ of the boundary of a vertex $V$ of $G$.
By Lemma 3.2, $p(\alpha )$ intersects $p(K^-)$.
Moreover, since the loop $\gamma \cup p(\alpha )$ bounds a disk $E$ in $S$, $|p(\alpha )\cap p(K^-)|=1$ and $\gamma$ meets exactly one of $e_1$ and $e_2$, say $e_1$.
Thus a loop $l=\partial N(\partial E; E)-\partial E$ intersects $\tilde{K}$ in two points.
Since $\tilde{K}$ is prime, $int E$ intersects $p(K)$ in an embedded arc.
Then, there are two posibilities for $e_1$, $e_1\subset f$ or $e_1\subset V$.
In the formar case, $f\cup E$ bounds a pair of a 3-ball and an unknotted arc, and an isotopy of $F$ along the pair eliminates $\alpha$.
In the later case, $f\cup E$ bounds a 3-ball , and an isotopy of $F$ along the 3-ball eliminates $\alpha$.
These contradict the minimality of $|(F\cap D^+)-(e_1\cup e_2)|$.
\end{proof}

Hence we have a condition that $G$ has at least two vertices, every vertex has valency at least two, and all faces of $G$ in $F$ are disks.
Next, we pay attention to a face of $G$ in $F$.

\begin{lemma}
For any face $f$, the cycle $\partial F$ can not be oriented.
\end{lemma}

\begin{proof}
If all corners of $f$ do not meet $e_1\cup e_2$, then this is same to Lemma 2.4.

If exactly one corner of $f$ meets $e_1$ or $e_2$ at one point, then $f$ and some $K^+_i$ have the intersection number $\pm 1$,
or a vertex which meets $f$ along the corner intersects some $K^-_k$ in one point.
Since $p(\partial f)$ and $p(K^-)\cap p(K^+_i)$ must have the intersection number zero, $\partial f$ is bounded by a loop of $G$ consisting of a vertex and an edge $\alpha$, and $p(\alpha )$ intersects $p(K^-)$ in one point.
Then Lemma 3.4 gives the conclusion.

If some corners of $f$ meet both $e_1$ and $e_2$, then the corners of $f$ have the intersection number zero with $p(K)$ because $F$ and $K$ have the intersection number zero.
In such a situation, we have a contradiction same as the proof of Lemma 2.4.
\end{proof}

By Lemma 3.5, starting a face $f$ of $G$ in $F$ whose closure is a disk, we can get a loop $l$ in $F-K$ with $|lk(l,K)|\ge 2$.
But this is imposible because any loop in $F-K$ is null-homotopic in $S^3-K$ or has linking number $\pm 1$ with $K$.
This finishes the proof of Theorem 1.4.




\end{document}